\documentstyle[amssymb,12pt,fleqn]{article}
\def\la{\langle}\def\ra{\rangle}

\def\nos#1{\leavevmode\unskip #1\ignorespaces }
\def\ns{\nos{-}}

\title{A Note on the Solvablity of Groups
\thanks{Project supported by the National Natural Science Foundation of
China.}
\author{Shiheng Li and Wujie Shi\\  \small (School of Mathematic
Science, Suzhou University, Suzhou 215006, Peoples Republic of China)\\
}}
\date{}

\begin{document}
\maketitle

\begin{abstract}Let $M$ be a maximal
subgroup of a finite group $G$ and $K/L$ be a chief factor such that
$L\leq M$ while $K\nsubseteq M$. We call the group $M\cap K/L$ a
$c$\ns section of $M$. And we define $Sec(M)$ to be the abstract
group that is isomorphic to a $c$\ns section of $M$. For every
maximal subgroup $M$ of $G$, assume that $Sec(M)$ is supersolvable.
Then any composition factor of $G$ is isomorphic to $L_2(p)$ or
$Z_q$, where $p$ and $q$ are primes, and $p\equiv\pm 1(mod\ 8)$.
This result answer a question posed by ref. \cite{WL}.
\end{abstract}


In \cite{WL} the following question was posed:

{\large\bf Question.} For every maximal subgroup $M$ of a group $G$
assume that $Sec(M)$ is supersolvable. Is $G$ solvable?

For solving this Question we need some definitions and lemmas.

{\bf Definition 1}\cite[Definition 1.1]{WL}. Let $M$ be a maximal
subgroup of a finite group $G$ and $K/L$ be a chief factor such that
$L\leq M$ while $K\nsubseteq M$. We call the group $M\cap K/L$ a
$c$\ns section of $M$.

We say that there is a unique class of subgroups $U$ in a group
$G$, if every subgroup isomorphic to $U$ is conjugate to $U$ in
$G$.

{\bf Lemma 1}\cite[Lemma 1.1]{WL}. For any maximal subgroup $M$ of
a group $G$, there is a unique $c$\ns section of $M$ up to
isomorphism.

By Lemma 1, it is reasonable to introduce the following:

{\bf Definition 2}\cite[Definition1.2(1)]{WL}. Given a maximal
subgroup $M$ of a group $G$, we define $Sec(M)$ as the abstract
group that is isomorphic to a $c$\ns section of $M$.

Similar to the proof of \cite[${\mathrm{II}}$, 5.3(a)]{H} and of
(1) in the proof of \cite[${\mathrm{II}}$, 5.3(b)]{H}, we get
Lemma 2.

{\bf Lemma 2.} (a) For $n\neq 4$, $1<k<n$, $A_n$ has no subgroup
of index $k$.

(b) For every subgroup $U$ of index $n$ in $A_n$ , there exists an
automorphism $\alpha$ of $A_n$ such that $U^\alpha=V_1$, where
$V_1=\{g\in A_n|1^g=1\}\cong A_{n-1}.$

{\bf Lemma 3.} If $n\neq 6$, the subgroups of index $n$ in
alternative group $A_n$are conjugate in $A_n$.

{\bf Proof.} Since $n\neq 6$, by \cite[${\mathrm{II}}$,5.5(a) and
5.3(b)]{H} there exist $n$ subgroups of index $n$ in $S_n$ and they
are conjugate in $S_n$. Let $U$ is a subgroup in $A_n$ of index $n$.
Since $Aut(A_n)\cong S_n\ (n\neq 6)$ and $A_n\ (n\neq 4)$ is simple,
$N_{S_n}(U)\cong S_{n-1}$ is a subgroup of $S_n$ of index $n$ by
Lemma 2(b). Moreover $S_{n-1}$ has only one subgroup of index 2.
Thus there exactly exist $n$ subgroups of index $n$ in $A_n$. Then
these subgroups are conjugate in $A_n$.\\

The following example gives a negative answer for the above
Question.

{\bf Example} Groups $$G:=PGL_2(p)$$ are the counterexamples of
the Question, where $p\in\mathbb{P}$ and $p\equiv \pm1(mod 8)$.

{\bf Proof}\ \ In the following, $p$ is always a prime with
$p\equiv \pm 1(mod\ 8)$. Let $G=PGL_2(p)$, $K=L_2(p)$. We know
that $G=K.\la \alpha\ra$, where $\alpha$ is an outer automorphism
of order 2 of $K$.

By \cite[II,8.16]{H}, $K$ has two conjugacy classes of elementary
abelian subgroups $T$ of type $(2,2)$ and $|N_K(T)|=24$. By
\cite[II,8.27]{H}, $PGL_2(p)$ is contained in $L_2(p^2)=:F$, then
$N_K(T)\leq N_G(T)\leq N_F(T)$. Thus, $|N_G(T)|=24$ since
$|N_F(T)|=24$ by \cite[II,8.16]{H} again. It follows that
$|G:N_G(T)|=2|K:N_K(T)|$, then all the elementary abelian
subgroups of type $(2,2)$ in $K$ are conjugate in $G$. Therefore,
$\alpha\in G\backslash K$ interchanges the two conjugacy classes
of elementary abelian subgroups of type $(2,2)$ in $L_2(p)$.
Especially, $T^\alpha$ and $T$ can't be conjugate in $K$.

Evidently, $G\vartriangleright K\vartriangleright 1$ is the unique
chief series of $G$. Let $M$ be a maximal subgroup of $G$ (denoted
$M<\cdot\ G$). If $M=K$, then $Sec(M)=1$. So we may assume $K\nleq
M$. Then $Sec(M)=K\cap M$ and $G=KM$. Hence, there exist elements
$k\in K$ and $m\in M$ such that $\alpha=km$. Then $T^m$ and $T$
cannot be conjugate in $K$.

By \cite[${\mathrm{II}}$,8.27]{H}, the maximal subgroups of $K$
are isomorphic to $S_4$, $A_5$, or supersolvable. Hence, if
$Sec(M)$ is not supersolvable, then $Sec(M)$ is isomorphic to
$S_4$, $A_5$, or $A_4$. Assume that $T\leq Sec(M)$, then $\la
T,T^m\ra\leq K\cap M=Sec(M)$. On the other hand, all elementary
abelian subgroups $T$ of type $(2,2)$ in $A_5$, $S_4$, or $A_4$
are conjugate in them respectively. Then $T^m$ and $T$ are
conjugate in $K$, a contradiction.

Therefore, $Sec(M)$ is supersolvable for each $M<\cdot\ G$.\\

Now, we give a complete answer to the Question:

{\large\bf Theorem.} Let $G$ be a finite group. For every maximal
subgroup $M$ of $G$, assume that $Sec(M)$ is supersolvable.Then
any composition factor of $G$ is isomorphic to $L_2(p)$ or $Z_q$,
where $p$ and $q$ are primes, and $p\equiv\pm 1(mod\ 8)$.\\

We need the following lemma.

{\bf Lemma 4.} Suppose that
$G=SL(n,p^f), f>1$ and $P\in Syl_p(G)$. Let $L=L_n(p)$ and $P_1\in
Syl_p(L)$. Then $N_G(P)$ and
$N_L(P_1)$ are non-supersolvable.\\

{\bf Proof}. By \cite[II,7.1]{H}, we may assume that $P$ consists
of the following matrices:

$$\left(\begin{array}{cccccc}
1&0&0&\cdots&0&0\\ a_{21}&1&0&\cdots&0&0\\
a_{31}&a_{32}&1&\cdots&0&0\\
\cdots&\cdots&\cdots&\cdots&\cdots&\cdots\\
a_{n-1,1}&a_{n-1,2}&a_{n-1,3}&\cdots&1&0\\
a_{n,1}&a_{n,2}&a_{n,3}&\cdots&a_{n,n-1}&1
\end{array}\right).$$

Then the normalizer of $P$ in $G$ is consists of the following
matrices:

$$\left(\begin{array}{cccccc}
a_{11}&0&0&\cdots&0&0\\ a_{21}&a_{22}&0&\cdots&0&0\\
a_{31}&a_{32}&a_{33}&\cdots&0&0\\
\cdots&\cdots&\cdots&\cdots&\cdots&\cdots\\
a_{n-1,1}&a_{n-1,2}&a_{n-1,3}&\cdots&a_{n-1,n-1}&0\\
a_{n,1}&a_{n,2}&a_{n,3}&\cdots&a_{n,n-1}&a_{n,n}
\end{array}\right),$$ where $a_{11}a_{22}\cdots a_{n,n}=1.$

Let $E=diag\{1,1,\cdots,1\}$ and $D\in N_G(P)$. Let $E_{ij}$ be a
matrix having a lone 1 as its $(i,j)$-entry and all other entries
0, where $i\neq j$. Then
\begin{eqnarray*}& &(D)^{-1}(E+aE_{n,1})D\\
&=&E+aD^{-1}E_{n,1}D\\
&=&E+a\left(\begin{array}{cccccc}
a_{11}&0&0&\cdots&0&0\\ a_{21}&a_{22}&0&\cdots&0&0\\
a_{31}&a_{32}&a_{33}&\cdots&0&0\\
\cdots&\cdots&\cdots&\cdots&\cdots&\cdots\\
a_{n-1,1}&a_{n-1,2}&a_{n-1,3}&\cdots&a_{n-1,n-1}&0\\
a_{n,1}&a_{n,2}&a_{n,3}&\cdots&a_{n,n-1}&a_{n,n}
\end{array}\right)^{-1}E_{n,1}\\
& &\left(\begin{array}{cccccc}
a_{11}&0&0&\cdots&0&0\\ a_{21}&a_{22}&0&\cdots&0&0\\
a_{31}&a_{32}&a_{33}&\cdots&0&0\\
\cdots&\cdots&\cdots&\cdots&\cdots&\cdots\\
a_{n-1,1}&a_{n-1,2}&a_{n-1,3}&\cdots&a_{n-1,n-1}&0\\
a_{n,1}&a_{n,2}&a_{n,3}&\cdots&a_{n,n-1}&a_{n,n}
\end{array}\right)\\
&=&E+a\left(\begin{array}{ccccc}
a_{11}^{-1}&0&\cdots&0&0\\ a'_{21}&a_{22}^{-1}&\cdots&0&0\\
\cdots&\cdots&\cdots&\cdots&\cdots\\
a'_{n-1,1}&a'_{n-1,2}&\cdots&a_{n-1,n-1}^{-1}&0\\
a'_{n,1}&a'_{n,2}&\cdots&a'_{n,n-1}&a_{n,n}^{-1}
\end{array}\right)\\
& &E_{n,1}\left(\begin{array}{ccccc}
a_{11}&0&\cdots&0&0\\ a_{21}&a_{22}&\cdots&0&0\\
\cdots&\cdots&\cdots&\cdots&\cdots\\
a_{n-1,1}&a_{n-1,2}&\cdots&a_{n-1,n-1}&0\\
a_{n,1}&a_{n,2}&\cdots&a_{n,n-1}&a_{n,n}
\end{array}\right)\\
&=&E+aa_{n,n}^{-1}E_{n,1}\left(\begin{array}{ccccc}
a_{11}&0&\cdots&0&0\\ a_{21}&a_{22}&\cdots&0&0\\
\cdots&\cdots&\cdots&\cdots&\cdots\\
a_{n-1,1}&a_{n-1,2}&\cdots&a_{n-1,n-1}&0\\
a_{n,1}&a_{n,2}&\cdots&a_{n,n-1}&a_{n,n}
\end{array}\right)\\
&=&E+aa_{n,n}^{-1}a_{11}E_{n,1}.
\end{eqnarray*}

And it is evident that $N=\{E+aE_{n,1}|a\in GF(p^f)\}$ is a
subgroup of the group of all matrices of $n$ rank with respect to
matrix multiplication. Hence $N$ is a normal subgroup of $N_G(P)$.

In case $n>2$, we can get that $a_{n,n}=1$. Thus, if $a\neq 0$ then
$aa_{n,n}^{-1}a_{11}=aa_{11}$ runs over $GF(p^f)$ when $a_{11}$ runs
over $GF(p^f)$. Therefore $N$ is a minimal normal subgroup of
$N_G(P)$. If $n=2$, then $a_{n,n}^{-1}=a_{11}$. Thus
$aa_{n,n}^{-1}a_{11}=aa_{11}^2$. It follows that
$|\{aa_{11}^2|a_{11}\in GF(p^f)\}|\geq p^f-1$ if $p=2$ and
$|\{aa_{11}^2|a_{11}\in GF(p^f)\}|=\frac{p^f-1}{2}$ if $p$ is odd
for $a\neq 0$. Then $N$ is also a minimal normal subgroup of
$N_G(P)$ when $n=2$ since $|N|=p^f$. Therefore, $N_G(P)$ is not
supersolvable for $f>1$.

Since $L\cong G/Z(G)$ and $P\cong P_1$, similar to above, we also
have that $N_L(P_1)$ is not supersolvable for $f>1$.\\

{\bf Proof of the Theorem}. Suppose that the Theorem is false. Then
there exists a minimal counterexample. Let $G$ be a minimal
counterexample.

(1)\ \ $G$ has the unique minimal normal subgroup $N$ and $N\cong
N_1\times N_2 \times\cdots\times N_t$, where $N_1\cong\cdots\cong
N_t$ is a non-abelian simple group.

By Lemma 1, $Sec(M/N)$ is supersolvable for every normal subgroup
$M/N$ of $G/N$, thus $G/N$ satisfies the hypotheses of the Theorem
. Hence any composition factor of $G/N$ is isomorphic to $L_2(p)$
or $Z_q$. Thus $G$ has the unique minimal normal subgroup $N$
since any composition factor of $G/L\cap K$ is isomorphic to a
composition factor of $G/K$ or $G/L$, where $K$ and $L$ are normal
subgroup of $G$. And $N$ is non-abelian since $G$ is a minimal
counterexample. Therefore, $N\cong N_1\times N_2\times\cdots\times
N_t$, where $N_1\cong\cdots\cong N_t$ is non-abelian simple group.

(2)\ \ Without loss of generality, we may assume that $G\cong
Aut(N_1)$, $N=Soc(G).$

By \cite[18.14]{DH}, $G$ might be considered as a subgroup $Inn(G)$
(the inner automorphism group of $N$) of $Aut(N)\cong
Aut(N_1)\wr_{nat}S_n$. Let $M_1$ be a maximal subgroup of $Aut(N_1)$
such that $N_1\nleq M_1$, then $Aut(N)=N(M_1\wr_{nat}S_n)$ and
$N\nleq M_1\wr_{nat}S_n$. Thus $G=G\cap
Aut(N)=N(G\cap(M_1\wr_{nat}S_n))$ and $G$ has a maximal subgroup $M$
that contains $G\cap(M_1\wr_{nat}S_n)$. Hence $G=MN$ and $N\nleq M$.
Then $Sec(M)=M\cap N\geq G\cap(M_1\wr_{nat}S_n)\cap N\geq N_1\cap
M_1$. By the hypothesis of the theorem, $Sec(M)$ is supersolvable
and thus $N_1\cap M_1$ is supersolvable. Hence $Aut(N_1)$ satisfies
the hypothesis of the theorem. On the other hand, since
$Aut(N_1)/N_1$ is solvable by Schreier conjecture and $G/N$
satisfies the Theorem, any composition of $G$ is isomorphic to
$L_2(p)$, $Z_q$, or $N_1$, where $p,q\in\mathbb{P}$ and $p\equiv\pm
1(mod\ 8)$.

Hence, if $G$ is a counterexample of the Theorem then $Aut(N_1)$
is also a counterexample of the Theorem; if $G$ isn't a
counterexample of the Theorem then $Aut(N_1)$ isn't a
counterexample of the Theorem. Therefore, we may assume that
$G\cong Aut(N_1)$ and $N=N_1$.

In the following, $G\cong Aut(N_1)$, $N=Soc(G).$

(3) Every maximal subgroup $M$ of $G$ such that $N\nleq M$ is
soluble. And $N\cap M=Sec(M)$ is supersolvable.

Since $M/Sec(M)\cong MN/N=G/N$ is soluble and $Sec(M)$ is
supersolvable by the hypothesis of the Theorem, then $M$ is
soluble.

Evidently, by (3) we get (4):

(4) If there is a proper subgroup $H$ of $G$ such that $HN=G$,
then $H$ is soluble and $H\cap N$ is supersolvable.

(5)\ \ Let $1\neq U$ be a proper subgroup of $N$. If there is a
unique class of subgroups $U$ in $N$ then $U$ is supersolvable.

Set $$\Gamma:=\{U^n|n\in N\}.$$
Since there is a unique class of
subgroups $U$ in $N$, both $G$ and $N$ transitively act on $\Gamma$
by conjugation. Therefore, $G=N_G(U)N$ by Frattini Argument. Then,
by (4) $N_G(U)$ is soluble and $U\leq N_G(U)\cap
N$ is supersolvable.\\

We analyse it case by case:

Case(A) $N=A_n,\ n\geq 5$. If $n\neq 6$ then $G=S_n$. Evidently,
both $G$ and $N$ act on $\{1,2,\cdots,n\}$ transitively, by Frattini
argument, $G=NS_{n-1}$. Thus, by (4), $A_{n-1}=N\cap S_{n-1}$ is
supersolvable, a contradiction. If $n=6$, from \cite{CNPW}, $G$
doesn't satisfy the hypothesis of the Theorem by (3).

Case(B) $N=L_2(q), q=p^f>3$. If $f>1$, then $N_G(P)$ is not
supersolvable by Lemma 4. If $p^2\not\equiv 1(mod\ 16)$,i.e.
$p\not\equiv\pm1(mod\ 8)$, then there is a unique class of
subgroups $A_4$ in $L_2(p)$. Thus $N\ncong L_2(q),f>1$, or $f=1$
and $p\not\equiv\pm1(mod\ 8)$ by (5). On the other hand, by
Example, $G$ satisfies the Theorem if $N\cong L_2(p),
p\equiv\pm1(mod\ 8)$.

Case(C) $N=L_n(q), n>2$. Since $L_3(2)\cong L_2(7)$, we suppose
$(n,q)\neq (3,2)$. By \cite[13.2]{A},
$G=N_G(M)PGL_n(q)=N_G(M)L_n(q)$, where we choose $M\cong
PGL_1(q)\times PGL_{n-1}(q)$. Evidently $M\cap N$ has a section
isomorphic to $L_{n-1}(q)$. But by (4) $M\cap N\leq N_G(M)\cap N$ is
supersolvable. Then $(n,q)=(3,2)$.

Case(D) $N=U_n(q), n\geq 3$. Since $G\leq P\Gamma L_n(q)$, both $G$
and $U_n(q)$ transitively act on the set of nonsingular subspaces of
dimension $i$ by Witt's theorem. Define $N_i$ to be the stabilizer
of a nonsingular space of dimension $i$ in $G$. Then $G=N_1N$ by
Frattini argument and thus $N_1$ is soluble. On the other hand, both
$N_1$ and $N_1\cap N$ have a section isomorphic to $U_{n-1}(q)$.
Then $(n,q)=(3,2)$ by (4). If $N=U_3(2)$, $G$ does not satisfy the
hypothesis by \cite{CNPW}.

Case(E) $N=PSp_4(q),q=2^f>2$. Let $P\in Syl_2(N)$. From \cite[\S
14]{A}, there is a $M<\cdot G$ such that $M\cap L=N_N(P)$. By (3),
$N_N(P)$ is supersolvable and then $N_N(P)=P\times H$, where $H$ is
a 2-complement of $N_N(P)$. From \cite[5.1.7(b)]{LPS}, we get that
$|N_N(P)|=q^4(q-1)^2$. So $|H|=(q-1)^2$. Let $r$ is the largest
prime divisor of $|H|$ and $R\in Syl_r(H)$. Then $PR=P\times R$ is a
nilpotent Hall $\{2,r\}$-subgroup of $N$. By \cite[II,9.24,b)]{H},
we can consider $PSp_2(q^2)$ as a subgroup of $N$. Let
$T=PSp_2(q^2)$ and $P_1\in Syl_2(T),\ R_1\in Syl_r(T)$. Since
$|N|=q^4(q+1)^2(q-1)^2(q^2+1)$ and $|T|=q^2(q-1)(q+1)(q^2+1)$, from
\cite{Wie}, we may assume that $P_1R_1\leq PR$. Thus $C_T(P_1)\geq
R_1>1$, contrary to $C_T(P_1)=1$.

Case(F) $N=PSp_{2m}(q),m>2$ or $m=2$ and $q$ odd. Then $G\leq
P\Gamma L_n(q)$. Both $G$ and $U_n(q)$ transitively act on the set
of totally singular $i$-subspaces for each $i$ by Witt's theorem.
Define $P_i$ to be the stabilizer of a totally singular $i$-space in
$G$. Then $G=P_1N$ by Frattini argument. On the other hand, both
$P_1$ and $P_1\cap N$ have a section isomorphic to
$PSp_{2(m-1)}(q)$. Then $PSp_{2(m-1)}(q)$ is soluble by (4), a
contradiction.

Case(G) $N=P\Omega_{2m+1}(q),m\geq 3,q$ odd. Define $N_i$ to be the
stabilizer of a nonsingular space of dimension $i$ in $G$. Similar
to Case $N=U_n(q), n\geq 3$, we get that $N_1$ is soluble and $N_1$
has a section isomorphic to $P\Omega_{2m}^+(q)$ or
$P\Omega_{2m}^-(q)$, a contradiction.

Case(H) $N=P\Omega_{8}^+(q)$. By \cite[Proposition 2.2.1]{K}, there
is a unique class of subgroups $R_1$ in $N$, where $R_1$ to be the
stabilizer of a totally singular 1-space in $N$. Then $R_1$ is
supersolvable by (5) but $R_1$ has a section isomorphic to
$P\Omega_{6}^+(q)$, a contradiction.

Case(I) $N=P\Omega_{2m}^+(q),m>4$. Then $G\leq P\Gamma L_n(q)$. Both
$G$ and $N$ transitively act on the set of totally singular
$i$-subspaces for each $i$ by Witt's theorem. Define $P_i$ to be the
stabilizer of a totally singular $i$-space in $G$. Then $G=P_2N$ by
Frattini argument. On the other hand, $P_2$ has a section isomorphic
to $P\Omega_{2(m-1)}^+(q)$. Then $PSp_{2(m-1)}(q)$ is soluble by
(4), a contradiction.

Case(J) $N=P\Omega_{2m}^-(q),m\geq 4$. Similar to Case
$N=P\Omega_{2m}^+(q),m>4$.

Case(K) $N$ is an exceptional group of Lie type:

Subcase(a) $N={}^2B_2(q),q=2^{2m+1}$. Let $P\in Syl_p(N)$. From
\cite{Su}, we know that $N_N(P)$ is a Frobenius group of order
$2^{2(2m+1)}(2^{2m+1}-1)$ and $N_N(P)$ is supersolvable by (5), a
contradiction.

Subcase(b) $N={}^2G_2(q),q=3^{2m+1}$. In this case the Sylow
2-subgroup $P$ of $N$ is abelian. In addition, $N_N(P)$ is
supersolvable by (5). Thus $N_N(p)=C_N(P)$. Hence $N$ is
2-nilpotent by the well-known theorem of Burnside and then $N$ is
soluble by the odd order theorem, a contradiction.

Subcase(c) $N=G_2(q)$. From \cite[Table 1]{LS}, there is a unique
class of subgroups $SL_3(q)$ in $G_2(q)$ if $3\nmid q$; there is a
unique class of subgroups ${}^2G_2(q)$ in $G_2(q)$ if $q=3^{2m+1}$;
there is a unique class of subgroups $G_2(q^m)$ in $G_2(q)$ if
$q=3^{2m}$, which contradicts (5).

Subcase(d) $N={}^3D_4(q)$. From \cite[Table 1]{LS}, there is a
unique class of subgroups $G_2(q)$ in $N$, which contradicts (5).

Subcase(e) $N={}^2F_4(q), q=2^{2m+1}>2,$ or ${}^2F_4(2)'$. If
$N={}^2F_4(2)'$, from \cite{CNPW}, we get that $G$ does not satisfy
the hypothesis by (3). If $N={}^2F_4(q), q=2^{2m+1}>2$, by
\cite[Proposition 2.12]{M}, there is just one class of subgroups
$L_2(25)$ in $N$. Hence $L_2(25)$ is supersolvable by (5), a
contradiction.

Subcase(f) $N=F_4(q)$. If $q$ is odd, from \cite[Table 1]{LS}, there
is a unique conjugacy class of subgroups $B_4(q)$ in $G_2(q)$. If
$q=2^m$, by \cite[2.8,2.9,and 2.11]{Gu}, $N_N(S)=SH$ is a group of
order $q^{24}(q-1)^4$ and $C_H(S)=1$, where $S\in Syl_2(N)$ and
$|H|=(q-1)^4$. However, by (5), $N_N(S)$ is supersolvable and thus
2-nilpotent. Then $N_N(S)=S\times H$, contrary to $C_H(S)=1$ when
$q>2$. When $q=2$, from \cite{CNPW}, we get that $G$ can't satisfy
the hypotheses by (3).

Subcase(g) $N={}^2E_6(q)$, $E_6(q)$,$E_7(q)$ or $E_8(q)$. In this
case, there is a unique class of non-soluble subgroups in $N$ by
\cite[Proposition 3.1, 4.1,5.1,6.1]{LS}, which contradicts (5).

Case(L) $N$ is a sporadic simple groups.

If the outer automorphism group $Out(N)=1$, then $G=N$ cannot
satisfy the hypothesis by \cite[VI, 9.6]{H}. Hence, we only consider
the sporadic simple groups $N$ such that $Out(N)=2$ since $|G/N|\leq
2$ and then, from \cite{CNPW}, we get that $G$ does not satisfy the
hypothesis by (3).

Therefore the minimal counterexample is not existence and the
theorem is proved.\\

{\bf Acknowledgment}. The authors are grateful to Profs. Shirong Li
and Xianhua Li for their helps. The authors are also grateful to the
referee for his suggestions.

\end{document}